\documentclass[a4paper, 10pt]{article}

\usepackage{amsmath}
\usepackage{amsthm}
\usepackage{amssymb}
\usepackage{stmaryrd}

\usepackage[vietnamese,USenglish]{babel}
\usepackage[utf8]{inputenc}

\usepackage{csquotes} 

\usepackage{enumerate}
\usepackage[hidelinks]{hyperref}

\usepackage{cleveref}

\theoremstyle{definition}
\newtheorem{definition}{Definition}
\theoremstyle{plain}
\newtheorem{theorem}[definition]{Theorem}
\newtheorem{corollary}[definition]{Corollary}

\newcommand{\N}{\mathbb{N}}
\newcommand{\R}{\mathbb{R}}
\newcommand{\Card}{\mathrm{Card}}
\newcommand{\Ker}{\mathrm{Ker}}
\renewcommand{\Im}{\mathrm{Im}}

\begin{document}

\title{When lattice cubes meet affine subspaces:\\a short note}
\author{\foreignlanguage{vietnamese}{\textsc{Nguyễn} Lê Thành Dũng}\\
  LIPN, CNRS \& Université Paris 13\\
\texttt{nltd@nguyentito.eu}}
\date{}
\maketitle

In this note, we give simple proofs of what seem to be folklore results:

\begin{theorem}
  \label{thm}
  Let $S \subseteq \R^d$ ($d \in \N$) be an affine subspace and $n \in \N$. The
  intersection of $S$ with the $d$-dimensional lattice cube of side length $n$
  -- that is, $S \cap \{0,\ldots,n-1\}^d$ -- has cardinality at most $n^{\dim
    S}$.
\end{theorem}
\begin{corollary}
  \label{cor}
  The minimum number of $k$-dimensional affine subspaces of $\R^d$ necessary to
  cover the lattice cube $\{0,\ldots,n-1\}^d$ is $n^{d-k}$.
\end{corollary}

The proof of \Cref{thm} is essentially an undergraduate-level linear algebra
exercise. That said, it is not immediately obvious.

\paragraph{Context}

The analogous question to \Cref{cor} for \emph{linear} subspaces (a.k.a.\ vector
subspaces) is raised as an open problem in the
book~\cite{DBLP:books/daglib/0017422} (to be precise, Problem~6 in
Section~10.2). At the time of writing, the state of the art on this problem is
the recent paper~\cite{DBLP:journals/dcg/BalkoCV19}.

Although both~\cite{DBLP:books/daglib/0017422,DBLP:journals/dcg/BalkoCV19}
contain many references to papers on related problems, we have not managed to
find a proof for the affine case in the literature. Presumably it has been
deemed too trivial by the authors of~\cite{DBLP:books/daglib/0017422} to deserve
inclusion, since they merely write \enquote{Covering by linear subspaces instead
  of affine ones is more difficult}. We suspect that a proof might have been
written in~\cite{Talata}, but that paper does not seem to be accessible on the
Internet.

The motivation for writing this note arose because of some recent research by
N.\ K.\ Blanchard and S.\ Kachanovich
\cite{Blanchard2019Countingauthorisedpaths} making use of (a weaker bound than)
\Cref{thm}.

\paragraph{Proofs}

\emph{Notations:} We write $(e_1, \ldots, e_d)$ for the canonical basis of
$\R^d$. We denote the lattice cube by $C(n,d) = \{0,\ldots,n-1\}^d$.

\begin{proof}[Proof of \Cref{cor} from \Cref{thm}]
  Let $S_1,\ldots,S_m$ be $k$-dimensional affine subspaces covering $C(n,d)$.
  The theorem gives us $\Card(S_i \cap C(n,d)) \leq n^k$. By summing over $i =
  1, \ldots, m$, we get $n^d = \Card(C(n,d)) \leq mn^k$, hence $n^{d-k} \leq m$.
  Conversely, the bound can be reached by the naive covering which uses the
  subspaces $(\R e_1 \oplus \ldots \oplus \R e_k) + a_{k+1} e_{k+1} + \ldots +
  a_d e_d$ for each of the $n^{d-k}$ choices $(a_{k+1},\ldots,a_d) \in
  \{0,\ldots,n-1\}^{d-k}$.
\end{proof}

\begin{proof}[Proof of \Cref{thm}]
  Let $k = \dim S$ ($0 \leq k \leq d$). Let $D$ be the \emph{direction} of $S$
  (i.e. the unique linear subspace of $\R^d$ that can be obtained by translating
  $S$). Let us choose any linear map $p : \R^d \to E$ (for some vector space
  $E$) such that $\Ker(p) = D$, e.g.\ the projection onto the orthogonal
  complement of $D$.

  Since the image of a spanning set by a linear map spans its range, the set
  $\{p(e_1), \ldots, p(e_d)\}$ spans $\Im(p)$. Therefore, it can be reduced to a
  basis of $\Im(p)$, which we may take without loss of generality to be $p(e_1),
  \ldots, p(e_{d-k})$ (the number of vectors in this basis is indeed
  $\dim(\Im(p)) = \dim(\R^d) - \dim(D) = d-k$ by the rank-nullity theorem).

  Now, this means that $p$ induces a linear isomorphism from $\R e_1 \oplus
  \ldots \oplus \R e_{d-k}$ to $\Im(p)$, since it maps a basis of the former to
  a basis of the latter. In particular, $p$ is injective when restricted to $B =
  \{x_1 e_1 + \ldots + x_{d-k} e_{d-k} | x_i \in \{0,\ldots,n-1\} \} = C(n,d)
  \cap (\R e_1 \oplus \ldots \oplus \R e_{d-k})$.

  Let $v \in \R^d$. By linearity, $p$ is also injective on any $B + v$. This
  entails that $\Card((B + v) \cap S) = \Card(p((B+v) \cap S)) \leq
  \Card(p(S))$. And since $D = \Ker(p)$, $p(D) = \{0\}$ so $p(S)$ is a singleton
  (again by linearity, since $S = D + u$ for any $u \in S$). Putting this
  together, $(B + v) \cap S$ is either empty or a singleton.

  To conclude, observe that $C(n,d)$ is the union of $B + a_1 e_{d-k+1} + \ldots
  + a_k e_d$ for $a_i \in \{0,\ldots,n-1\}$. Since the intersection of each of
  these $n^k$ sets with $S$ has at most one point, we obtain the desired bound.
\end{proof}

\bibliographystyle{alpha}
\bibliography{lattice-subspace}

\end{document}